\documentclass[graybox]{svmult}
\usepackage{mathptmx}       
\usepackage{helvet}         
\usepackage{courier}        
\usepackage{type1cm}        
\usepackage{makeidx}        
\usepackage{graphicx}        
\usepackage{multicol}        
\usepackage[bottom]{footmisc} 
\usepackage{mathrsfs}
\usepackage{amsmath}
\usepackage{amsbsy}
\usepackage{amstext}

\usepackage{amssymb}

\begin{document}

\title{Recovering Composition Algebras from 3D Geometric Algebras}

\author{Daniele Corradetti}

\institute{Daniele Corradetti \at Universidade do Algarve, \at 
Campus de Gambelas \at \email{a55944@ualg.pt}
\and Grupo Fisica Matematica \at Instituto Superior Tecnico, Lisboa \at \email{d.corradetti@gmail.com}}  
\maketitle

\abstract{Generalized Hurwitz theorem states that there are fifteen composition
algebras: seven unital, six para-unital, and two non-unital algebras.
In this article we explore the recovery of such algebras from 3D Geometric
Algebras. Different involutions, such as reversion, inversion, Clifford
conjugation, and full grade inversion, are introduced in order to
recover the norm of all composition algebras. A special attention
is given to composition algebras of dimension 8, i. e. octonions,
para-octonions and Okubo algebra, for which the introduction of a
different product is needed.}
\section{Introduction}

A composition algebra $A$ is an algebra for which the bilinear product
is compatible the norm, i.e., a
\begin{align}
n\left(x\cdot y\right) & =n\left(x\right)n\left(y\right),\label{eq:comp(Def)-1-1}
\end{align}
for every $x,y\in A$. Composition algebras, which include real numbers
$\mathbb{R},$ , complex numbers $\mathbb{C}$, quaternions $\mathbb{H}$,
and octonions $\mathbb{O}$, are of particular interest due to their
unique properties and the way they can model physical phenomena. Composition
algebras are completely classified and can be divided in \emph{unital},
i.e. with $1\in A$ such that $x\cdot1=1\cdot x=x$; \emph{para-unital},
i.e. with an involution $x\longrightarrow\overline{x}$ and a paraunit
$\boldsymbol{1}\in A$ such that $x\cdot\boldsymbol{1}=\boldsymbol{1}\cdot x=\overline{x}$,
and, finally, \emph{non-unital}, i.e. intended here for algebras that
do not possess neither a unit element nor a paraunit element. 

As a consequence of the Generalized Hurwitz Theorem \cite{ElDuque Comp},
only fifteen composition algebras exist for any given field: seven
are unital and called \emph{Hurwitz algebras}; seven are paraunital
and termed \emph{para-Hurwitz algebras}; finally, there are two non-unital
and 8-dimensional, known as the \emph{Okubo algebras} $\mathcal{O}$
and $\mathcal{O}_{s}$

\medskip{}

\begin{center}
\begin{tabular}{|c|c|c|c|c|c|c|c|c|c|c|c|c|}
\hline 
\textbf{Hurwitz}  & \textbf{O.}  & \textbf{C.}  & \textbf{A.}  & \textbf{Alt.}  & \textbf{F.}  &  & \textbf{p-Hurwitz}  & \textbf{O.}  & \textbf{C.}  & \textbf{A.}  & \textbf{Alt.}  & \textbf{F.}\tabularnewline
\hline 
\hline 
$\mathbb{R}$  & Yes  & Yes  & Yes  & Yes  & Yes  &  & $p\mathbb{R}\cong\mathbb{R}$  & Yes  & Yes  & Yes  & Yes  & Yes\tabularnewline
\hline 
$\mathbb{C}$, $\mathbb{C}_{s}$  & No  & Yes  & Yes  & Yes  & Yes  &  & $p\mathbb{C}$, $p\mathbb{C}_{s}$  & No  & Yes  & No  & No  & Yes\tabularnewline
\hline 
$\mathbb{H}$,$\mathbb{H}_{s}$  & No  & No  & Yes  & Yes  & Yes  &  & $p\mathbb{H}$,$p\mathbb{H}_{s}$  & No  & No  & No  & No  & Yes\tabularnewline
\hline 
$\mathbb{O}$,$\mathbb{O}_{s}$  & No  & No  & No  & Yes  & Yes  &  & $p\mathbb{O}$,$p\mathbb{O}_{s}$  & No  & No  & No  & No  & Yes\tabularnewline
\hline 
\end{tabular}
\par\end{center}

\begin{center}
\begin{tabular}{|c|c|c|c|c|c|}
\hline 
\textbf{Okubo}  & \textbf{O.}  & \textbf{C.}  & \textbf{A.}  & \textbf{Alt.}  & \textbf{F.}\tabularnewline
\hline 
$\mathcal{O}$,$\mathbb{\mathcal{O}}_{s}$  & No  & No  & No  & No  & Yes\tabularnewline
\hline 
\hline
\end{tabular}
\medskip{}

\par\end{center}
Table 1.: \emph{This table summarize the properties of different composition algebras: totally ordered (O.), commutative (C.), associative (A.), alternative (Alt.), flexible (F.).}
\section{Three 8-dimensional composition algebras}

Let $A$ be unital composition algebra. Given a norm one has a \emph{polarization}
of it given by 
\[
\left\langle x,y\right\rangle =n\left(x+y\right)-n\left(x\right)-n\left(y\right),
\]
and then given the presence of the unit one has an involution $x\longrightarrow\overline{x}$
called \emph{conjugation} given by $\overline{x}=\left\langle x,1\right\rangle 1-x$.
The conjugation then satisfy the notable relation $x\cdot\overline{x}=n\left(x\right)1.$
In unital algebras the definition of a norm and the definition of
the conjugation is a strictly related process and, furthermore, it
is uniquely defined in unital composition algebras. Once the conjugation
over an unital composition algebra is established, then one can define
a new product
\[
x\bullet y=\overline{x}\cdot\overline{y},
\]
and since $n\left(x\right)=n\left(\overline{x}\right)$, one has that
$n\left(x\bullet y\right)=n\left(\overline{x}\cdot\overline{y}\right)=n\left(x\right)n\left(y\right),$
and thus the algebra $\left(A,\bullet,n\right)$ is again a composition
algebra. Nevertheless, since $\ x\bullet 1 =1\bullet x=\overline{x}$, the algebra is not unital anymore, but para-unital, i.e. $\left(A,\bullet,n\right)$ is the para-Hurwitz algebra of $A$. The same process can be achieved using an order three automorphisms $\tau$, i.e. $\tau^{3}=\text{id}$,
then defining a new product $*$ such that 
\[
x*y=\tau\left(\overline{x}\right)\cdot\tau^{2}\left(\overline{x}\right),
\]
and for which $\left(A,*,n\right)$ is again composition algebra which, for a suitable $\tau$, is non-unital \cite{Pe69,Ok78,KMRT,ElDuque Comp}.
From these considerations one can obtain the three 8-dimensional division
composition algebras that are deeply linked one another, but none
of them is isomorphic to the others, i.e. \medskip{}

\begin{center}
\begin{tabular}{|c|c|c|c|}
\hline 
Property  & $\mathbb{O}$  & $p\mathbb{O}$  & $\mathcal{O}$\tabularnewline
\hline 
\hline 
Unital  & Yes  & No  & No\tabularnewline
\hline 
Paraunital  & No  & Yes  & No\tabularnewline
\hline 
Alternative  & Yes  & No  & No\tabularnewline
\hline 
Flexible  & Yes  & Yes  & Yes\tabularnewline
\hline 
Composition  & Yes  & Yes  & Yes\tabularnewline
\hline 
Automorphism  & $\text{G}_{2}$  & $\text{G}_{2}$  & $\text{SU}_{3}$\tabularnewline
\hline 
\end{tabular}\medskip{}
\par\end{center}

Notably, all 8-dimensional composition algebras are obtainable from one another with a deformation of the bilinear product involving an order two antihomomorphism (conjugation) and a suitable order three automorphism $\tau$.
\medskip{}
\begin{center}
\begin{tabular}{|c|c|c|c|}
\hline 
Algebra  & $\left(\mathcal{O},*,n\right)$  & $\left(p\mathbb{O},\bullet,n\right)$  & $\left(\mathbb{O},\cdot,n\right)$\tabularnewline
\hline 
\hline 
$x*y$  & $x*y$  & $\tau\left(x\right)\bullet\tau^{2}\left(y\right)$  & $\tau\left(\overline{x}\right)\cdot\tau^{2}\left(\overline{y}\right)$\tabularnewline
\hline 
$x\bullet y$  & $\tau^{2}\left(x\right)*\tau\left(y\right)$  & $x\bullet y$  & $\overline{x}\cdot\overline{y}$\tabularnewline
\hline 
$x\cdot y$  & $\left(e*x\right)*\left(y*e\right)$  & $\left(\boldsymbol{1}\bullet x\right)\bullet\left(y\bullet\boldsymbol{1}\right)$  & $x\cdot y$\tabularnewline
\hline 
\end{tabular}
\par\end{center}
\medskip{}
Table 2.: \emph{Given the same 8-dimensional vector space, all three 8-dimensional composition algebras are obtainable one another through the definition of a new product as above.In all three cases the vector space and the norm remains the same}
\medskip{}
Similar relations hold for the split companion of these algebras,
i.e., $\mathcal{O}_{s}$, $p\mathbb{O}_{s}$, $\mathbb{O}_{s}$ (see
\cite{CMZ24b,CMZ24,ElDuque Comp}). Nonetheless, it is vital to note
that while transitioning from one algebra to another is feasible,
these algebras are not isomorphic. For example, while the octonions
$\mathbb{O}$ are alternative and unital, para-octonions $p\mathbb{O}$
are nor alternative nor unital but do have a para-unit. In contrast,
the Okubo algebra $\mathcal{O}$ is non-alternative and only contains
idempotent elements.
\section{Complex and Quaternionic Algebras from 3D Geometric algebras }
Let $\mathcal{G}\left(p,q\right)$ be the the Clifford algebra $Cl_{p,q}\left(\mathbb{R}^{3}\right)$,
i.e., where the quadratic form $Q$ has signature $\left(p,q\right)$.
Let $\left\{ e_{1},e_{2},e_{3}\right\} $ be an orthonormal basis
for $\mathbb{R}^{3}$, then $\mathscr{B}=\left\{ 1,e_{1},e_{2},e_{3},e_{1}e_{2},e_{2}e_{3},e_{1}e_{3},e_{1}e_{2}e_{3}\right\} $
is a basis for $\mathcal{G}\left(p,q\right)$. Without loss of generality
let assume that in such basis $Q$ is diagonal, i.e., $Q=\text{diag\ensuremath{\left(\lambda_{1},\lambda_{2},\lambda_{3}\right)}}$
with $\lambda_{1},\lambda_{2},\lambda_{3}\in\left\{ \pm1\right\} $.
Then one has the following relations 
\begin{align}
e_{i}^{2} & =\lambda_{i},\left(e_{i}e_{j}\right)^{2}=-\lambda_{i}\lambda_{j},\left(e_{1}e_{2}e_{3}\right)^{2}=-\lambda_{1}\lambda_{2}\lambda_{3},\label{eq:SquaredRelations}
\end{align}
where we intended $i\neq j$. Moreover, one also verifies that

\begin{equation}
\begin{array}{ccc}
\left(e_{1}e_{2}\right)\left(e_{2}e_{3}\right) & = & \lambda_{2}\left(e_{1}e_{3}\right),\\
\left(e_{2}e_{3}\right)\left(e_{1}e_{3}\right) & = & \lambda_{3}\left(e_{1}e_{2}\right),\\
\left(e_{1}e_{3}\right)\left(e_{1}e_{2}\right) & = & \lambda_{1}\left(e_{2}e_{3}\right).
\end{array}\label{eq:QuaternionicRelations}
\end{equation}

Let $Rot\left(\mathcal{G}\right)$ be the (split)-quaternionic subalgebra
of rotors with basis $\mathscr{R}=\left\{ 1,e_{1}e_{2},e_{2}e_{3},e_{1}e_{3}\right\} $
and $Ps\left(\mathcal{G}\right)$ the (split)-complex subalgebra with
basis $\left\{ 1,e_{1}e_{2}e_{3}\right\} .$ Then, one has a correspondence
between $\mathcal{G}\left(p,q\right)$ and biquaternionic algebra
$Ps\left(\mathcal{G}\right)\otimes Rot\left(\mathcal{G}\right)$ summarized
in the following table, i.e.,

\medskip{}

\begin{center}
\begin{tabular}{|c|c|c|c|}
\hline 
Algebra  & Tensor product  & $Rot\left(\mathcal{G}\right)$  & $Ps\left(\mathcal{G}\right)$\tabularnewline
\hline 
\hline 
$\mathcal{G}\left(3,0\right)$  & $\mathbb{C}\otimes\mathbb{H}$  & $\mathbb{H}$  & $\mathbb{C}$\tabularnewline
\hline 
$\mathcal{G}\left(2,1\right)$  & $\mathbb{C}_{s}\otimes\mathbb{H}_{s}$  & $\mathbb{H}_{s}$  & $\mathbb{C}_{s}$\tabularnewline
\hline 
$\mathcal{G}\left(1,2\right)$  & $\mathbb{C}\otimes\mathbb{H}_{s}$  & $\mathbb{H}_{s}$  & $\mathbb{C}$\tabularnewline
\hline 
$\mathcal{G}\left(0,3\right)$  & $\mathbb{C}_{s}\otimes\mathbb{H}$  & $\mathbb{H}$  & $\mathbb{C}_{s}$\tabularnewline
\hline 
\end{tabular}
\par\end{center}

\medskip{}
\noindent where the product of the algebra for the tensor product is intended
$\left(z_{1}\otimes q_{1}\right)\left(z_{2}\otimes q_{2}\right)=z_{1}z_{2}\otimes q_{1}q_{2}$
for every $z_{1},z_{2}\in C$ and $q_{1},q_{2}\in C$ with $C\in\left\{ \mathbb{C},\mathbb{C}_{s}\right\} $
and $Q\in\left\{\mathbb{H},\mathbb{H}_{s}\right\} $. 
For many practical purposes it is vital to define a specific norm, or, in
an equivalent manner, an involution. A classical example is that of
integral system over the algebra, which are usually defined as 
\[
n\left(x\right),tr\left(x\right)\in\mathbb{Z},
\]
where $tr\left(x\right)=\left\langle x,1\right\rangle $ (see for
example \cite{Di23}). With this in mind,let $c$ be the pseudoscalar
$c=e_{1}e_{2}e_{3},$ and 
\begin{equation}
c_{+}=\frac{1+c}{2},\,\,\,\,\,c_{-}=\frac{1-c}{2},
\end{equation}
so that $x=x_{+}+x_{-},$ where $x_{+}=xc_{+}$ and $x_{-}=xc_{-}$.
We define axiomatically the following involutions: \emph{reversion},
\emph{inversion}, \emph{Clifford conjugation} and \emph{full grade
inversion}. Correspondence between involutions of 3D geometric algebra
and different conjugation types defined over biquaternionic algebras
is defined as follow

\medskip{}

\begin{center}
\begin{tabular}{|c|c|c|c|}
\hline 
\textbf{Name}  & \textbf{Involution}  & \textbf{Definition}  & \textbf{Biquaternionic}\tabularnewline
\hline 
\hline 
Reversion  & $x\longrightarrow x^{\dagger}$  & $\left\langle x\right\rangle _{1}^{\dagger}=\left\langle x\right\rangle _{1},\left(xy\right)^{\dagger}=y^{\dagger}x^{\dagger}$  & $x^{\dagger}=\overline{z}\otimes\widetilde{q}$\tabularnewline
\hline 
Inversion  & $x\longrightarrow\overline{x}$  & $\overline{x}=x_{+}-x_{-}$  & $\overline{x}=\overline{z}\otimes q$\tabularnewline
\hline 
Clifford conjugation  & $x\longrightarrow\widetilde{x}$  & $\widetilde{x}=\overline{x}^{\dagger}$  & $\widetilde{x}=z\otimes\widetilde{q}$\tabularnewline
\hline 
Full grade inversion  & $x\longrightarrow x^{*}$  & $x^{*}=x_{+}^{\dagger}-x_{-}.$  & only Octonionic\tabularnewline
\hline 
\end{tabular}
\par\end{center}

\medskip{}
It is worth noting that different involutions give rise to different
norms (and thus to different integral systems). Indeed, defining the
norm $\widetilde{n}\left(x\right)=x\widetilde{x}\in\mathbb{C}$ then
$\left(\mathcal{G}\left(3,0\right),\widetilde{n}\right)$ is a composition
algebra isomorphic to the quaternionic algebra over the complex field,
i.e., $\mathbb{H}_{\mathbb{C}}$, with integral unit elements forming
a $D_{4}\subset\mathbb{R}^{4}$ root systems; while $n_{\dagger}\left(x\right)=xx^{\dagger}\in\mathbb{R}$
then $\left(\mathcal{G}\left(3,0\right),n_{\dagger}\right)$ is not
a composition algebra with integral unit elements forming a $D_{4}\oplus D_{4}\subset\mathbb{R}^{8}$
root system which can be then projected to an $F_{4}$ root system
in $\mathbb{R}^{4}$. 

As for the Hurwitz and para-Hurwitz complex and quaternionic algebra,
then one has the following:
\begin{align}
\left(Rot\left(\mathcal{G}\left(p,q\right)\right),\widetilde{n}\right) & \cong\left(\mathbb{H},\cdot,n\right),\left(p,q\right)\in\left\{ \left(3,0\right),\left(0,3\right)\right\} ,\\
\left(Rot\left(\mathcal{G}\left(p,q\right)\right),\widetilde{n}\right) & \cong\left(\mathbb{H}_{s},\cdot,n\right),\left(p,q\right)\in\left\{ \left(2,1\right),\left(1,2\right)\right\} ,\\
\left(Ps\left(\mathcal{G}\left(p,q\right)\right),\overline{n}\right) & \cong\left(\mathbb{C},\cdot,n\right),\left(p,q\right)\in\left\{ \left(3,0\right),\left(1,2\right)\right\} ,\\
\left(Ps\left(\mathcal{G}\left(p,q\right)\right),\overline{n}\right) & \cong\left(\mathbb{C}_{s},\cdot,n\right),\left(p,q\right)\in\left\{ \left(0,3\right),\left(2,1\right)\right\} ,
\end{align}
while the para-Hurwitz algebras $p\mathbb{H},p\mathbb{H}_{s}$ and
$p\mathbb{C},p\mathbb{C}_{s}$ are achieved substituting the Clifford
product with the following two products 
\begin{align}
x \widetilde{\bullet} y & =\widetilde{x}\widetilde{y},\\
x\overline{\bullet}y & =\overline{x}\overline{y},
\end{align}
 to the previous subalgebras $Rot\left(\mathcal{G}\right)$ and $Ps\left(\mathcal{G}\right)$.

\section{8-Dimensional Composition Algebras from 3D Geometric algebras }

All 8-Dimensional composition algebras are not associative so one
has to define a new product $\bullet$ and $\bullet_{-}$ which was inspired by \cite{Hi22}. Let $\tau$
be the order three map defined by
\begin{equation}
\tau:\begin{cases}
\begin{array}{cccc}
1\longrightarrow & 1, & e_{1}e_{2}\longrightarrow & e_{1}e_{2},\\
e_{1}\longrightarrow & -\frac{1}{2}\left(e_{1}-\sqrt{3}e_{3}\right), & e_{2}e_{3}\longrightarrow & -\frac{1}{2}\left(e_{2}e_{3}-\sqrt{3}e_{2}\right),\\
e_{2}\longrightarrow & -\frac{1}{2}\left(e_{2}-\sqrt{3}e_{2}e_{3}\right), & e_{1}e_{3}\longrightarrow & e_{1}e_{3}\\
e_{3}\longrightarrow & -\frac{1}{2}\left(e_{3}-\sqrt{3}e_{1}\right), & e_{1}e_{2}e_{3}\longrightarrow & e_{1}e_{2}e_{3}
\end{array}\end{cases}
\end{equation}

Then over $\mathcal{G}\left(p,q\right)$ we define the following products

\begin{align}
x\cdot y & =x_{+}y_{+}+\widetilde{y_{-}}x_{-}+y_{-}x_{+}+x_{-}\widetilde{y_{+}},\\
x\bullet y & =\widetilde{x_{+}}\widetilde{y_{+}}+\widetilde{y_{-}}x_{-}-y_{-}\widetilde{x_{+}}-x_{-}y_{+},\\
x*y & =\tau\left(\widetilde{x_{+}}\right)\tau^{2}\left(\widetilde{y_{+}}\right)+\tau^{2}\left(\widetilde{y_{-}}\right)\tau\left(x_{-}\right)-\tau^{2}\left(y_{-}\right)\tau\left(\widetilde{x_{+}}\right)-\tau\left(x_{-}\right)\tau^{2}\left(y_{+}\right),
\end{align}
which give rise to the Octonionic, para-Octonionic and Okubo algebra.
Moreover, switching the sign to the second and third term of the above
product, one obtains other three products, i.e. $\cdot_{-},\bullet_{-}$
and $*_{-}$, for which one recovers the split version of the previous
algebras. Thus, one has the following \medskip{}

\begin{center}
\begin{tabular}{|c|c|c|c|c|c|c|}
\hline 
Algebra  & $\left(\mathcal{G}\left(p,q\right),\cdot\right)$  & $\left(\mathcal{G}\left(p,q\right),\bullet\right)$  & $\left(\mathcal{G}\left(p,q\right),*\right)$  & $\left(\mathcal{G}\left(p,q\right),\cdot_{-}\right)$  & $\left(\mathcal{G}\left(p,q\right),\bullet_{-}\right)$  & $\left(\mathcal{G}\left(p,q\right),*_{-}\right)$\tabularnewline
\hline 
\hline 
$\mathcal{G}\left(3,0\right)$  & $\mathbb{O}$  & $p\mathbb{O}$  & $\mathcal{O}$  & $\mathbb{O}$  & $p\mathbb{O}$  & $\mathcal{O}$\tabularnewline
\hline 
$\mathcal{G}\left(2,1\right)$  & $\mathbb{O}_{s}$  & $p\mathbb{O}_{s}$  & $\mathcal{O}_{s}$  & $\mathbb{O}_{s}$  & $p\mathbb{O}_{s}$  & $\mathcal{O}_{s}$\tabularnewline
\hline 
$\mathcal{G}\left(1,2\right)$  & $\mathbb{O}_{s}$  & $p\mathbb{O}_{s}$  & $\mathcal{O}_{s}$  & $\mathbb{O}_{s}$  & $p\mathbb{O}_{s}$  & $\mathcal{O}_{s}$\tabularnewline
\hline 
$\mathcal{G}\left(0,3\right)$  & $\mathbb{O}_{s}$  & $p\mathbb{O}_{s}$  & $\mathcal{O}_{s}$  & $\mathbb{O}_{s}$  & $p\mathbb{O}_{s}$  & $\mathcal{O}_{s}$\tabularnewline
\hline 
\end{tabular} \medskip{}
\par\end{center}

\begin{flushleft}
To obtain the norm of the above composition algebras, one has to use
the full grade inversion, $n^{*}\left(x\right)=xx^{*}$, for which
the full expression is 
\begin{align}
n^{*}\left(x\right) & =x_{0}^{2}+x_{1}^{2}\left(\lambda_{1}\lambda_{2}\right)+x_{2}^{2}\left(\lambda_{2}\lambda_{3}\right)+\\
 & x_{3}^{2}\left(\lambda_{1}\lambda_{3}\right)+x_{4}^{2}\lambda_{1}+x_{5}^{2}\lambda_{2}+x_{6}^{2}\lambda_{3}+x_{7}^{2}\lambda_{1}\lambda_{2}\lambda_{3},
\end{align}
 for any $x=x_{0}+x_{1}e_{1}e_{2}+x_{2}e_{2}e_{3}+x_{3}e_{1}e_{3}+x_{4}e_{1}+x_{5}e_{2}+x_{6}e_{3}+x_{7}e_{1}e_{2}e_{3}$.
Finally, one has that the integral units of $\left(\mathcal{G}\left(3,0\right),n^{*}\right)$
is an $E_{8}$ root system as it is well known from the study of integral
octonions \cite{Co46}. 
\par\end{flushleft}

\section{Acknowledgments}
The author wishes to express sincere gratitude to Alessio Marrani, Francesco Zucconi, Richard Clawson and Klee Irwin for their invaluable discussions and insights that significantly contributed to the development of this work. Special thanks are also extended to Eckhard Hitzer and Carlos Castro Perelman for their help in partecipating to ICHAA 2024 and the stimulating conversations.



\begin{thebibliography}{CMZ24b}
\bibitem[CMZ24]{CMZ24}Corradetti D., Marrani A., Zucconi F., \emph{A
minimal and not alternative realisation of the Cayley plane}, Ann.
Univ. Ferrara (2024) doi.org/10.1007/s11565-024-00498-5

\bibitem[MCZ23]{MCZ23-1}Marrani A., Corradetti D., Zucconi F., \emph{Physics
with non-unital algebras? An invitation to the Okubo algebra} (2023)
arXiv:2309.17435.

\bibitem[CRI23]{CRI23}Corradetti D., Clawson R., Irwin K., \emph{All
Hurwitz Algebras from 3D Geometric Algebras} (2023) arXiv:2311.02269

\bibitem[CMZ24b]{CMZ24b}Corradetti D., Marrani A. \& Zucconi F.,
\emph{Collineation groups of octonionic and split-octonionic planes
}arXiv:2311.11907

\bibitem[Co46]{Co46}Coxeter, H. S. M., \emph{Integral Cayley numbers},
Duke Math. J., 13 (4) (1946), 561--578

\bibitem[Di23]{Di23}Dickson L. \emph{Algebras and their arithmetics},
Chicago, Ill.: The University of Chicago Press 1923.

\bibitem[El18]{ElDuque Comp}Elduque A., \emph{Composition algebras};
in \emph{Algebra and Applications I: Non-associative Algebras and
Categories}, Chapter 2, pp. 27-57, edited by Abdenacer Makhlouf, Sciences-Mathematics,
ISTE-Wiley, London 2021.

\bibitem[Hi22]{Hi22}Hitzer, E. \emph{Extending Lasenby's embedding
of octonions in space-time algebra $Cl\left(1,3\right)$, to all three-
and four dimensional Clifford geometric algebras}. Math Meth Appl
Sci. 2022; 1-24. doi:10.1002/mma.8577

\bibitem[Hur]{Hurwitz98}Hurwitz A., \emph{Uber die Komposition der
quadratischen Formen von beliebig vielen Variablen}, Nachr. Ges. Wiss.
Gottingen, 1898.

\bibitem[KMRT]{KMRT}Knus M.-A., Merkurjev A., Rost M. and Tignol
J.-P., \emph{The book of involutions.American Mathematical Society
Colloquium Publications} 44, American Mathematical Society, Providence,
RI, 1998.

\bibitem[Ok78]{Ok78}Okubo S., \emph{Deformation of Pseudo-quaternion
and Pseudo-octonion Algebras}, Hadronic J. 1 (1978) 1383.

\bibitem[Pe69]{Pe69}Petersson H.P. , \emph{Eine Identitat funften
Grades, der gewisse Isotope von Kompositions- Algebren genugen}, Math.
Z. 109 (1969), 217--23

\end{thebibliography}
\end{document}